\documentclass[a4paper,10pt]{article}
\usepackage{fullpage}
\usepackage{amsfonts} 
\usepackage{amssymb} 
\usepackage{latexsym} 
\usepackage{epsfig} 
\usepackage{amsmath}
\usepackage{amsthm}
\usepackage[all]{xy}
\usepackage{mathrsfs}

\newcommand{\zz}{\mathbb{Z}}
\newcommand{\qq}{\mathbb{Q}}
\newcommand{\cc}{\mathbb{C}} 
\newcommand{\pp}{\mathbb{P}} 
\newcommand{\oo}{\mathcal{O}} 
\newcommand{\mx}{\mathfrak{m}} 
\newcommand{\lra}{\longrightarrow}                         
\newcommand{\ra}{\rightarrow}
\newcommand{\xra}{\xrightarrow}

\newcommand{\cvf}[2][]{\frac{\partial #1}{\partial #2}} 
\newcommand{\gen}[1]{\left<#1\right>}%
\newcommand{\set}[1]{\left\{ #1 \right\}}%
\newcommand{\id}{\mathbb{I}}%

\newtheorem{theorem}{Theorem}[section]%
\newtheorem{lemma}[theorem]{Lemma}%
\newtheorem{corollary}[theorem]{Corollary}%
\newtheorem{proposition}[theorem]{Proposition}%
\theoremstyle{definition} 
\newtheorem{definition}[theorem]{Definition}%
\newtheorem{remark}[theorem]{Remark}%
\newtheorem*{acknowledgement}{Acknowledgement}

\DeclareMathOperator{\Lie}{Lie}%
\DeclareMathOperator{\Res}{Res}%
\begin{document}

\title{Some examples of non-massive Frobenius manifolds in Singularity Theory}%

\author{Ignacio de Gregorio}

%
%

%
%
\maketitle
\begin{abstract}
  Let $f,g:\cc^2\ra\cc$ be two quasi-homogeneous polynomials. We compute the
  $V$-filtration of the restriction of $f$ to any plane curve $C_t=g^{-1}(t)$ and show
  that the Gorenstein generator $dx\wedge dy/dg$ is a primitive form. Using results of
  A.~Douai and C.~Sabbah, we conclude that base space of the miniversal unfolding of
  $f_t:=f|_{C_t}$ is a Frobenius manifold. At the singular fibre $C_0$ we obtain a
  non-massive Frobenius manifold.
\end{abstract}

\section{Introduction}
\label{sec:introduction}

The axiomatic of Frobenius manifolds as originally defined by B.~Dubrovin as for example
in~\cite{Dub96}, represents the geometrisation of the celebrated WDVV or associativity
equations in topological quantum field theories (c.f.~\cite{DijVerVer91}). This
geometrisation made it plain evident that Frobenius manifold, and hence solutions to WDVV
equations, already existed in a very different branch of mathematics, namely singularity
theory and more particularly deformations of hypersurface singularities. This work had
been carried out by K.~Saito and M.~Saito nearly ten years before (see~\cite{KSai81},
\cite{KSai83} and \cite{MSai89}).

The other main source of Frobenius manifolds is quantum cohomology, where the solutions
are a priori, just formal series and can only be geometrised after some effort if at
all. A version of the mirror phenomenon is interpreted in this framework as an isomorphism
of two Frobenius manifolds, each coming from one of this two seemingly unrelated
sources. In this direction, we have the result of S.~Barannikov (\cite{Bar00})
establishing an isomorphism between the quantum cohomology of projective spaces and the
Frobenius manifold obtained by unfolding the function $x_0+\dots+x_n$ on the affine
variety $x_0\dots x_n=1$.

As the mirror of $\pp^n$ indicates, in order to find potential mirrors of algebraic
varieties, it is not enough to look at Frobenius manifolds produced by unfolding of germs
of isolated singularities. Global functions on affine varieties are needed. A.~Douai and
C.~Sabbah in \cite{DouSab03} have adapted the results of M.~Saito to this global affine
situation and, under some mild hypothesis, reduced the existence of Frobenius-Saito
structures on the base space of the miniversal unfolding to the existence of a primitive
form for the Gauss-Manin system. They used their results to exhibit Frobenius structures
for unfoldings of non-degenerate and convenient Laurent polynomials.

In this article, we construct Frobenius manifolds for unfoldings of quasi-homogeneous
functions on quasi-homogeneous plane curves. Let $f,g:\cc^2\ra\cc$ be quasi-homogeneous
polynomials with respect to the same weights. We regard $g$ as a family of plane curves
$C_t=g^{-1}(t)$, and consider the restriction $f_t:=f|_{C_t}$. We show that the
$V$-filtration of the Gauss-Manin system, and hence the spectral pairs, of $f_t$ can be
computed from $f_0$. In particular, the Gorenstein generator $\alpha:=dx\wedge dy/dg$
yields a primitive form with associated spectral number $0$. It follows from
\cite{DouSab03} that the base space of the miniversal deformation of $f_t$ can be endowed
with a Frobenius manifold structure. At $t=0$, the curve $C_0$ has an isolated
singularity. We can use the dualising module $\omega_{C_0}$ to define the Gauss-Manin
system of $f_0$ and the Grothendieck residue pairing to construct a non-massive Frobenius
manifold. 

The motivation behind the construction is the following remark: the unfolding of
$f=x^a+y^b$ on $C_t:xy=t$, $t\not=0$, is the mirror partner of the weighted projective
line $\pp(a,b)$ (for $a$ and $b$ coprimes, see~\cite{Mann05,CoaCorLeeTse06}). At $t=0$, the
multiplication and metric in our construction at the origin coincides with the orbifold
cohomology of $\pp(a,b)$.

\begin{acknowledgement}
  The author is indebted to C.~Sabbah for helpful discussions during the preparation of
  this article.
\end{acknowledgement}

\section{Preliminaries}
\label{sec:preliminaries}

Let us recall briefly how to obtain a Frobenius manifold from a meromorphic connection. We
closely follow C.~Sabbah (cf.~\cite{Sab98}).

Let $G\ra B$ be a vector bundle on manifold $B$ and let the rank of $G$ be equal to the
dimension of $B$, say $m$. Let $F$ denote the pull-back of $G$ via the projection
$\mathbb{P}^1\times B \ra B$. We further assume that $F$ is equipped with a flat
meromorphic connection $\widehat{\nabla}$ with logarithmic poles along $\set{0}\times B$
and poles of type $1$ along $\set{\infty}\times B$. From this initial data we obtain the
following objects:
\begin{list}{(\roman{enumi})}
  {\usecounter{enumi}
    \setlength{\leftmargin}{5\jot}}
\item the residual connection $\nabla$ on $G\ra B$: if $\tau$ denotes the coordinate on
  the affine chart $\pp^1 \setminus\set{\infty}$ the connection matrix for
  $\widehat{\nabla}$ is locally written as
  \begin{equation}
    \label{eq:1}
    \Omega^{\widehat{\nabla}}
    =\Omega_{\tau}\frac{d\tau}{\tau}+\sum_{i=1}^m\Omega_i du_i
  \end{equation}
  where $\Omega_\tau$ and $\Omega_i$ are matrices with holomorphic entries.  Here
  $(u_1,\dots,u_n)$ denotes a coordinate system on a neighbourhood in $B$. The residual
  connection $\nabla$ on $B$ is given by
  \begin{displaymath}
    \Omega^{\nabla}=\sum_{i=1}^m\Omega_i(0,u_1,\dots,u_n) du_i.
  \end{displaymath}
  and the integrability of $\widehat{\nabla}$ implies that of $\nabla$;
\item the residue endomorphism of $\widehat{\nabla}$, that is, an endomorphism $R_0$ of
  $F|_B$ given in local coordinates by $\Omega_{\tau}(0,u)$. The integrability of
  $\widehat{\nabla}$ implies that $R_0$ is covariantly constant with respect to $\nabla$,
  i.e., $\nabla R_0 = 0$;
\item an endomorphism $R_\infty$ of $F|_B$, defined (up to constant) by the choice of a
  coordinate $\theta$ in $\pp^1\setminus\set{0}$. Indeed, the connection at infinity has a
  pole of type $1$. If we use $\theta=\tau^{-1}$ as a coordinate in
  $\pp^1\setminus\set{0}$ we see from (\ref{eq:1}) that $\widehat{\nabla}$ is written near
  $\infty$ as
  \begin{displaymath}
    \frac{1}{\theta}\big(\Omega_\theta\frac{d\theta}{\theta}
    + \sum_{i=1}^m\Omega_i' du_i\big)
  \end{displaymath}
  where $\Omega_\theta=-\theta\Omega_\tau$ and $\theta^2\Omega_i'=\Omega_i$ have
  holomorphic entries. The matrix $\Omega_\theta(0,u_1,\dots,u_n)$ defines the
  endomorphism $R_\infty$ of $F|_B$. The coordinate $\theta$ (and hence $\tau$) will be
  kept fixed throughout this article.
\item The Higgs field $\Phi$, defined as follows. We decompose the connection
  $\widehat{\nabla}=\widehat{\nabla}'+\widehat{\nabla}''$ according to the decomposition
  of $1$-forms $\pi_{\pp^1\setminus\set{0}}^*\Omega^1_{\pp^1\setminus\set{0}} \oplus
  \pi_{B}^*\Omega^1_{B}$. We write $\widehat{\nabla}''=d_B + \Omega''$ and set
  $\Phi=(\theta\Omega'')|_{\theta=0}$. It also depends on the choice of the coordinate
  $\theta$ (up to constant).
\end{list}
The integrability of $\widehat{\nabla}$ implies the following relations between all of the
above objects:
\begin{align}
  \label{eq:2}
  \begin{split}
    &\nabla^2=0,~\nabla R_0 =0\\
    &\Phi\wedge\Phi=0,~[R_\infty,\Phi]=0\\
    &\nabla\Phi=0,~\nabla R_\infty + \Phi = [\Phi, R_0]
  \end{split}
\end{align}
Let $\mathcal{F}[*(\set{0}\times B)]$ denote the module of sections of $F$ with poles along
$\set{0}\times B$ and let $\mathbf{F}$ denote locally free $\oo_B[\theta]$-module
$(\pi_B)_*\mathcal{F}[*(\set{0}\times B]$. We further assume that $\mathbf{F}$ is equipped
with with a non-degenerate $\cc$-linear pairing
\begin{displaymath}
  S:\mathbf{F}\otimes\mathbf{F}\lra \theta\oo_B[\theta]
\end{displaymath}
satisfying
\begin{align}
  \label{eq:3}
  \begin{split}
    S(\theta m,m')&=\theta S(m,m')=S(m,-\theta m')\\
    \Lie_{\partial_\theta} S(m,m')&= S(\widehat{\nabla}_{\partial_\theta} m, m') +
    S(m,-\widehat{\nabla}_{\partial_\theta} m')\\
    \Lie_{\partial_t} S(m,m')&= S(\widehat{\nabla}_{\partial_t} m, m') +
    S(m,\widehat{\nabla}_{\partial_t} m')
  \end{split}
\end{align}
Expanding $S$ as a series in $\theta=0$ we get
\begin{displaymath}
  S(m,m')=\theta s_\infty^{1}(m,m') + \theta^2 s_\infty^{2}(m,m') + \dots.
\end{displaymath}
It can be checked that $s_\infty^1$ is a non-degenerate, symmetric pairing on
$\mathbf{F}/\theta\mathbf{F}$ which is metric with respect to the connection $\nabla$.
For a $\nabla$-horizontal section $\omega$ of $G$. We define its associate period mapping
$\varphi_\omega:TB\lra G$ by
\begin{displaymath}
  \varphi_\omega(\xi):=-\Phi(\xi)(\omega)
\end{displaymath}
We say that $\omega$ as above is primitive if
\begin{enumerate}
\item $\omega$ is an eigenvector of $R_0$ and
\item $\varphi_\omega$ is an isomorphism.
\end{enumerate}
If $\omega$ is a primitive form, we can define a $\oo_{B}$-algebra structure on $\Theta_B$
by setting
\begin{equation}
  \label{eq:4}
  \varphi_\omega(\xi\star\eta):=-\varphi^{-1}_\omega\Phi(\xi)\varphi_\omega(\eta)
\end{equation}
and we obtain:
\begin{theorem}
  \label{prop:1}
  (\cite{Sab98}) If $\omega$ is a primitive form, the triple $(B,\star,s_\infty)$ is a
  Frobenius manifold.
\end{theorem}
\begin{remark}
  \label{prop:2}
  We finish this section with a remark that simplifies enormously the construction of
  Frobenius manifolds from families of meromorphic connections.  Namely, if $B$ is simply
  connected, it is enough to check the existence of the primitive form at one single value
  of the parameter space $B$. This result is proved in a detailed manner in \cite{Sab98},
  but it goes back to the work of B.~Dubrovin on isomonodromic deformations.
\end{remark}

\section{Functions on curves}
\label{sec:deformations}

Let us recall the definition of the Milnor number of a function $f_0$ on a curve-germ given
by D.~Mond and D.~van~Straten~in \cite{MonStr01}.
\begin{definition}
  \label{prop:3}
  Let $(C,0)\hookrightarrow(\cc^{n},0)$ be a reduced curve-germ and let $f_0:(C,0)\ra
  (\cc,0)$ be a function non-constant on any branch. The Milnor number $\mu$ of $f_0$ is
  defined as
  \begin{equation}
    \label{eq:5}
    \mu:=\dim_\cc\frac{\omega_{C,0}}{\oo_{C,0}df_0}
  \end{equation}
  where
  $\omega_{C,0}=\text{Ext}^{n-1}_{\oo_{\cc^{n},0}}(\oo_{C,0},\Omega^{n}_{\cc^{n},0})$
  denotes the dualising module of $\oo_{C,0}$.
\end{definition}
\begin{remark}
  \label{prop:4}
  The authors in~\cite{MonStr01} show that if the curve is unobstructed (i.e. the second
  cotangent cohomology group $T^2_{C,0}$ vanishes) then the local Milnor numbers are
  preserved under flat deformation of $(C,0)$ and arbitrary deformation of $f_0$.
\end{remark}
In the case of complete intersection curves the Milnor number is relatively easy to
compute. If $(C,0)$ is a complete intersection curve defined by $g_1,\dots,g_n$, the
dualising module $\omega_{C,0}$ can be identified with the module of meromorphic $1$-forms
$\omega$ on $(C,0)$ such that $\omega\wedge dg_1\wedge\dots\wedge
dg_n\in\oo_{C,0}\otimes\Omega^{n+1}_{\cc^{n+1},0}$. It is therefore customary to write
$\omega_{C,0}=\oo_{C,0}\alpha$ where
\begin{equation}
  \label{eq:6}
  \alpha=\frac{dx_1\wedge\dots\wedge dx_{n+1}}{dg_1\wedge\dots\wedge dg_n}
\end{equation}
Given now $f_0:(C,0)\ra(\cc,0)$, let $f$ be a representative of $f_0$ in
$\oo_{\cc^n,0}$. We can write $df_0=J\alpha$ where $J$ is the Jacobian determinant of the
map $\varphi=(f,g_1,\dots,g_n):(\cc^{n+1},0)\ra(\cc^{n+1},0)$. Hence $\mu=\oo_{C,0}/(J)$
and using the L\^e-Greuel formula we see that
\begin{equation}
  \label{eq:7}
  \mu=\mu_{1}+\mu_{2}
\end{equation}
where $\mu_{1}$ denotes de Milnor number of $(C,0)$ and $\mu_{2}$ that of the
$0$-dimensional complete intersection defined by $\varphi$.

An unfolding of $f$ over $(B,0)=(\cc^m,0)$ is a function $F:(\cc^n\times B,0)\ra (\cc,0)$
together with fibration $(g,\id):(\cc^n\times B)\ra (\cc^{n-1}\times B,0)$ such that
$F|_{C_0}=f_0$. We say that $F$ is a {\em miniversal unfolding (resp. versal)} if the
Kodaira-Spencer map
\begin{equation}
  \label{eq:8}
  \Theta_{B,0}\ni\cvf{u_i}\mapsto
  \cvf[F]{u_i}\in\frac{\oo_{C_0\times B,0}}{(J)}
\end{equation}
is an isomorphism (resp. a surjection) of $\oo_{B,0}$-modules. Here $(u_1,\dots,u_m)$
denote coordinates on $(B,0)$. Notice that if $C_t$ is a Milnor fibre of an appropriate
representative of $g$, conservation of the Milnor number implies that the map
\begin{equation}
  \label{eq:9}
  \Theta_{B,0}\ni\cvf{u_i}\mapsto
  \cvf[F]{u_i}\in\frac{\oo_{C_t\times (B,0)}}{(J)}
\end{equation}
is also an isomorphism (resp. surjection). Hence the the restriction of $F$ to $C_t\times
(B,0)$$f_0$ is a miniversal deformation of $F|_{C_t}$ in the usual left-equivalence sense
for multigerms.

Notice also that the isomorphisms \eqref{eq:8} and \eqref{eq:9} induce structures of
$\oo_{B}$-algebras on the tangent sheaf $\Theta_B$. It is proved in \cite{deGre06} that
these multiplicative structures satisfy certain integrable condition turning them into
$F$-manifolds (see \cite{HerMan99,Her02}).

\subsection{The quasi-homogeneous case}
\label{sec:quas-funct-compl}

As noted in Remark~\ref{prop:4}, the Milnor number is locally preserved under
deformations. Here we wish to show that in the quasi-homogeneous case it is actually
globally preserved. Later, this will justify the use of algebraic forms to study the
Gauss-Manin system.

Most of the calculations that follow can be carried out for the case of complete intersections curve
singularities and we do so. However, our techniques can only be used to construct Frobenius
manifolds for functions on plane curves as it is in this case that we are able to extract
information about the spectrum of the restriction of the miniversal unfolding of $f_0$ to the Milnor
fibre of the singularity.  

Let us begin by introducing some notation that will be kept for the remainder of this article. Let
$\oo$ denote the polynomial ring $\cc[x_1,\dots,x_{n+1}]$. We make $\oo$ into a graded ring by
assigning the positive rational weight $p_i$ to the variable $x_i$. Homogeneity will always mean
homogeneity with respect to this grading. Let us be given
\begin{enumerate}
\item a polynomial map $g:\cc^{n+1}\ra\cc^n$ where $g_i$ is homogeneous of degree $e_i$,
  we denote the fibre over $t\in\cc^n$ by $C_t$ and suppose that the $0$-fibre $C_0$ is
  not smooth (see Remark~\ref{prop:5} below);
\item a homogeneous polynomial $f\in\oo$ of degree $1$, we write $f_t$ for the restriction
  of $f$ to the fibre $C_t$ and assume that $f_0$ is not constant on any branch of $C_0$.
\end{enumerate}
\begin{remark}
  \label{prop:5}
  The smooth case is exceptional as it is the only case for which $f$ belongs to its
  Jacobian algebra. On the other hand, the smooth case corresponds to the deformation of
  the $A_\mu$-singularity in one variable, and it is well-known that the base space of its
  miniversal deformation does have a Frobenius structure.
\end{remark}
Let $\alpha=dx_1\wedge\dots\wedge dx_{n+1}/dg_1\wedge\dots\wedge dg_n$ and let $\omega_g$
be the relative dualising module. As before, let $J$ be the Jacobian determinant of
$(f,g_1,\dots,g_n)$ so that $df=J\alpha$. As $J$ is also homogeneous, the only critical
point of $f_0$ is the origin and $\mu=\dim_\cc\oo/(g_1,\dots,g_n,J)$. The following
proposition shows that this is also the number of critical points of $f_t$. Let
$(t_1,\dots,t_n)$ be coordinates on the target space of $g$.
\begin{proposition}
  \label{prop:6}
  The $\cc[t_1,\dots,t_n]$-module $\oo/(J)$ is free of rank $\mu$.
\end{proposition}
\begin{proof}
  The module $\oo/(J)$ can be seen as a graded module over the graded ring
  $\cc[t_1,\dots,t_n]:=\cc[g_1,\dots,g_n]$. As
  $\omega_{C_0}/\oo_{C_0}df=\oo/(g_1,\dots,g_n,J)$ is a finite dimensional vector space it
  follows from the graded Nakayama lemma that $\oo/(J)$ is finitely generated (we recall
  that the graded version of Nakayama lemma does not require that the module $\oo/(J)$ be
  finitely generated). As $(g_1,\dots,g_n,J)$ is a regular sequence, the graded version of
  the Auslander-Buchsbaum formula tells us that $\oo/(J)$ is free as
  $\cc[t_1,\dots,t_n]$-module.
\end{proof}

\section{The Gauss-Manin system}
\label{sec:gauss-manin-system}

We keep the notation and hypothesis introduced in the previous section. We define the
(algebraic) Gauss-Manin system of $f$ relative to $g$ as the module
\begin{displaymath}
  \mathbf{G}:=\frac{\omega_g[\tau,\tau^{-1}]}{(d-\tau
    df\wedge)\oo[\tau,\tau^{-1}]}
\end{displaymath}
where $d$ denotes the relative differential with respect to $g$. It is a
$\cc[t_1,\dots,t_n,\tau,\tau^{-1}]$-module endowed with a partial integrable connection
with respect to $\partial_\tau$ defined as:
\begin{equation}
  \label{eq:10}
  \widehat{{\nabla}}_{\partial_\tau}[\omega]=[-f\omega]
\end{equation}
We also consider the {\em (relative) Brieskorn lattice $G$}, that is, the image of the
canonical map $\omega_g[\tau^{-1}]\ra\mathbf{G}$. It is a lattice of $\mathbf{G}$ as the
following proposition shows:
\begin{proposition}
  \label{prop:7}
  $G$ is a free $\cc[t_1,\dots,t_n,\tau^{-1}]$-module of rank $\mu$.
\end{proposition}
\begin{proof}
  According to Prop.~\ref{prop:6}, let $h_1,\dots,h_\mu$ be a basis of the free
  $\cc[t_1,\dots,t_n]$-module $\oo/J$ consisting of homogeneous elements. Let
  $\omega_i=h_i\alpha$ and let $\omega=a_0\alpha\in\omega_g$. Then there exist unique
  $c_1,\dots,c_\mu\in\cc[t_1,\dots,t_n]$ such that $a_0=c_1h_1+\dots+c_\mu h_\mu + a_0'J$,
  which implies that $\omega=c_1\omega_1+\dots+c_\mu\omega_\mu +
  a_0'df=c_1\omega_1+\dots+c_\mu\omega_\mu + \tau^{-1}da_0'$. Writing $da_0'=a_1\alpha$ we
  see that $\deg a_0 > \deg a_1$. The proposition follows by iteration.
\end{proof}

We begin by studying the action of $\partial_\tau$ on the holomorphic (algebraic) forms
$\Omega_g$. Recall that we are excluding the case in which $C_0$ is smooth.
\begin{lemma}
  \label{prop:8}
  Let $I=(g_1,\dots,g_n)$ and let $J_g$ be the ideal generated by all the maximal minors
  of the Jacobian matrix of $g$. The sequence
  \begin{equation}
    \label{eq:11}
    0\lra\frac{J_g+I}{I+(J)}
    \lra\frac{\oo}{I+(J)}\xra{f\cdot}\frac{\oo}{I+(J)}
    \lra\frac{\oo}{(f)+I+(J)}\lra 0
  \end{equation}
  is exact.
\end{lemma}
\begin{proof}
  Let $\mu_1$ be the Milnor number of $C_0$ and $\mu_2$ that of the $0$-dimensional
  complete intersection defined by $\varphi:=(f,g_1,\dots,g_n)$. We know that
  $\mu=\mu_1+\mu_2$. If $\widetilde{E}=\sum_{i=1}w_ix_i\cvf{x_i}$ denotes the Euler vector
  field of $\oo$, the homogeneity of $f$ and $g_i$ gives us
  $t\varphi(\widetilde{E})=f\cvf{s}+\sum_{i=1}^n e_i g_i \cvf{t_i}$. Applying Cramer's
  rule we obtain
  \begin{equation}
    \label{eq:12}
    w_i x_i J=(-1)^{i+1}f M_i \mod I
  \end{equation}
  where $M_i$ is the minor of the Jacobian matrix of $g$ obtained by deleting the $i$-th
  column.  From here we see that $fJ_g\subset I+(J)$. But we have an exact sequence
  \begin{equation}
    \label{eq:13}
    0\lra \frac{(f)+I+(J)}{(f)+I}\lra
    \frac{\oo}{(f)+I}\lra\frac{\oo}{(f)+I+(J)}\lra 0
  \end{equation}
  where the middle term has dimension $\mu_2+1$ (cf.~\cite{Loo84},~Prop.~5.12).  So that
  can also use (\ref{eq:12}) to conclude that
  \begin{equation}
    \label{eq:14}
    \dim_\cc \frac{\oo}{(f)+I+(J)}=\mu_2
  \end{equation}
  On the other hand (\ref{eq:12}) together with Nakayama lemma tells us that the first
  term of (\ref{eq:13}) has dimension $1$ so that (\ref{eq:14}) follows. Back to the
  original sequence (\ref{eq:11}) we conclude that the kernel of $f\cdot$ has dimension
  $\mu_2$. That $\mu_2$ is also the dimension of the first term of (\ref{eq:11}) follows
  from one more exact sequence:
  \begin{displaymath}
    0 \lra \frac{J_g+I}{I+(J)} \lra \frac{\oo}{I+(J)} \lra
    \frac{\oo}{J_g+I} \lra 0.
  \end{displaymath}
  The middle and last term of the above sequence have dimension $\mu=\mu_1+\mu_2$ and
  $\mu_1$ respectively (\cite{Loo84},~Prop.~9.10). Therefore the first term has dimension
  $\mu_2$ and the lemma follows.
\end{proof}
The following notation will be useful to describe the action of $\partial_\tau$ on
$\Omega_g$.

\vspace{\parsep}

\noindent{\bf Notation.}~We set $\mathbf{e}=\sum_{i=1}^{n} e_i$,
$\mathbf{p}=\sum_{i=1}^{n+1} p_i$ and for a homogeneous element $h\in \oo$, we define
\begin{displaymath}
  \nu(h):=\deg h + \mathbf{p} - \mathbf{e}.
\end{displaymath}
We will also write $\nu(\omega):=\nu(h)$ where $\omega=h\alpha$.
\begin{remark}
  \label{prop:9}
  Notice that for $\omega=h\alpha$ with $h$ homogeneous we have $\Lie_{\widetilde{E}}
  (\omega) = \nu(h)\omega$, where $\widetilde{E}$ denote the Euler vector field on
  $\oo$. Also, if $h\in J_g$ and as before we denote by $M_i$ the minor of the Jacobian
  matrix of $g$ obtained by deleting the $i$-th columns, then
  \begin{equation}
    \label{eq:15}
    \deg(h)\geq\min\set{\deg
      M_i:i=1,\dots,n+1}=\mathbf{e}-\mathbf{p}+\max\set{p_i:i=1,\dots,n+1}.
  \end{equation}
  It follows that $\nu(\omega)>0$.
\end{remark}
\begin{lemma}
  \label{prop:10}
  Let $\omega\in\Omega_g$ be a homogeneous $1$-form. Then, in $\mathbf{G}$ we have
  \begin{equation}
    \label{eq:16}
    \tau\partial_\tau [\omega]=-\nu(\omega)[\omega] + \sum_{j=1}^n t_j
    \omega_j + \tau\sum_{j=1}^n t_j \omega'_j
  \end{equation}
  with $\nu(\omega_i) \leq \nu(\omega) -e_j$ and $\nu(\omega'_j)\leq\nu(\omega) + 1 -
  e_j$.
\end{lemma}
\begin{proof}
  By linearity, we can assume that $\omega=hdx_{n+1}$ with $h$ homogeneous. As
  $dx_{n+1}=M_{n+1}\alpha$, we see that $\nu(\omega)=\deg h + p_{n+1}$. Let us introduce
  some helpful notation to carry out the calculation: $i_{\partial_{x_i}}$ denotes the
  contraction with respect to the vector field $\partial_{x_i}$,
  $\mathfrak{i}_{n+1}=i_{\partial_{n}}\circ\dots\circ i_{\partial_1}$ and
  $\mathfrak{i}_{n+1,j}
  =i_{\partial_{n}}\circ\dots\circ\widehat{i}_{\partial_j}\circ\dots\circ i_{\partial_1}$.
  Writing $dx_{n+1}=\mathfrak{i}_{n+1}V$ we have
  \begin{align*}
    -\tau\widehat{\nabla}_\tau [\omega] &= \tau\partial_\tau[h\mathfrak{i}_{n+1}V] =
    \tau[h\mathfrak{i}_{n+1}(fV)] = \tau[h\mathfrak{i}_{n+1}(df\wedge
    i_{\widetilde{E}}V)]\\
    &= \tau\sum_{j=1}^n(-1)^{j+1}\left[h (i_{\partial_j}df)\wedge
      \mathfrak{i}_{n+1,j}i_{\widetilde{E}}V\right]
    +(-1)^n\tau\left[hdf\wedge\mathfrak{i}_{n+1}i_{\widetilde{E}}V\right]\\
    &= \tau\sum_{j=1}^n (-1)^{j+1}\left[h(i_{\partial_j}df)\wedge
      \mathfrak{i}_{n+1,j}i_{\widetilde{E}}V\right]+\tau\left[hdf\wedge i_{\widetilde{E}}\mathfrak{i}_{n+1}V\right]\\
    &= [di_{\widetilde{E}}\omega]+\tau\sum_{j=1}^n
    (-1)^{j+1}\left[h(i_{\partial_j}df)\wedge
      \mathfrak{i}_{n+1,j}i_{\widetilde{E}}V\right]\\
    &= [\Lie_{\widetilde{E}}(\omega)-i_{\widetilde{E}}d\omega]+\tau\sum_{j=1}^n
    (-1)^{j+1}\left[h(i_{\partial_j}df)\wedge \mathfrak{i}_{n+1,j}i_{\widetilde{E}}V\right]\\
    &= \nu(\omega)[\omega] - [i_{\widetilde{E}}d\omega]+\tau\sum_{j=1}^n
    (-1)^{j+1}\left[h(i_{\partial_j}df)\wedge
      \mathfrak{i}_{n+1,j}i_{\widetilde{E}}V\right]
  \end{align*}
  Multiplying the second summand by $dg_1\wedge\dots\wedge dg_n$ and applying the
  determinant theorem we see that $i_{\widetilde{E}}d\omega=\sum_{j=1}^n t_j\omega_j$
  where $\nu(\omega_j)=\nu(\omega)-e_j$, and analogously for the other summand. This can
  also be seen by noticing that the expression (\ref{eq:16}) is homogeneous where
  $\deg\tau=-1$.
\end{proof}
\begin{corollary}
  \label{prop:11}
  For $\omega\in\Omega_{C_0}$, we have $\tau\partial_\tau[\omega]=-\nu(\omega)[\omega]$ in
  $\mathbf{G}_0=\mathbf{G}/\mx_{\cc^n,0}$.
\end{corollary}

\subsection{$V$-filtration and spectral numbers}
\label{sec:spectral-numbers}

For a fixed point $t\in\cc^n$, set $\mathbf{G}_t=\mathbf{G}/\mx_{\cc^n,t}$ and analogously
$G_t=G/\mx_{\cc^n,t}$. Let us recall the definition of the Malgrange-Kashiwara
$V_\bullet$-filtration for $\mathbf{G}_t$. It is the unique filtration
$V_\bullet(\mathbf{G}_t)$ indexed by $\qq$ such that
\begin{enumerate}
\item $V_\lambda(\mathbf{G}_t)$ is $\cc[\tau]$-free and
  $\cc[\tau,\tau^{-1}]\otimes_{\cc[\tau]} V_\lambda(\mathbf{G}_t)=\mathbf{G}_t$ for all
  $\lambda\in\qq$;
\item $\tau V_\lambda(\mathbf{G}_t) \subset V_{\lambda-1},~\partial_\tau
  V_\lambda(\mathbf{G}_t) \subset V_{\lambda+1}$ and
\item the action of $\tau\partial_\tau + \lambda$ is nilpotent on the quotient
  $\text{gr}^V_\lambda(\mathbf{G}_t) := V_\lambda(\mathbf{G}_t)/V_{<\lambda}(\mathbf{G}_t)$.
\end{enumerate}
Such a filtration exists and is unique (e.g.~\cite{Bjork93},~pg.~113). Moreover, there
exists a finite subset $A\subset[0,1)$ such that $\text{gr}^V_\lambda(\mathbf{G}_t) = 0$
for all $\lambda\not\in A+\zz$.

The filtration $V_\bullet(\mathbf{G}_t)$ induces a filtration on $G_t/\tau^{-1}G_t$. The
corresponding graded part is given by
\begin{displaymath}
  \text{gr}^V_\lambda(G_t/\tau^{-1}G_t):=\frac{V_\lambda(\mathbf{G}_t)\cap G_t}
  {V_\lambda(\mathbf{G}_t)\cap \tau^{-1}G_t + V_{<\lambda}(\mathbf{G}_t)\cap G_t}.
\end{displaymath}
Let $d(\lambda)$ denote the dimension as a complex vector space of
$\text{gr}^V_\lambda(G/\tau^{-1}G)$. The set of pairs $(\lambda,d(\lambda))$ for which
$d(\lambda)\not=0$ is called the {\em spectrum} of $(\mathbf{G}_t,G_t)$.

We can use lemmas~\ref{prop:8} and \ref{prop:10} to compute the $V_\bullet$-filtration of
the Gauss-Manin system of the function $f_0$ and for the case of plane curves, for any
$f_t$. The linear map
$(-f_0)\cdot:\omega_{C_0}/\oo_{C_0}df_0\ra\omega_{C_0}/\oo_{C_0}df_0$ is nilpotent and
homogeneous. Hence its Jordan basis induces a homogeneous basis of $G$ of the following
form:
\begin{align}
  \label{eq:17}
  \begin{split}
    &\left[\omega_1^{i}\right]=\left[(-f)^{i}\omega_1^{0}\right],~i=0,\dots,N_1\\
    &\left[\omega_2^{i}\right]=\left[(-f)^{i}\omega_2^{0}\right],~i=0,\dots,N_2\\
    &\hdots\\
    &\left[\omega_{M}^{i}\right]=\left[(-f)^{i}\omega_{M}^{0}\right],~i=0,\dots,N_{M}\\
    &\left[\omega_{M+1}^{0}\right],\dots,\left[\omega_{\mu_2}^{0}\right]
  \end{split}
\end{align}
It is helpful to set $\nu_i^j=\nu(\omega_i^j)$. Consider now the following change of basis
of $G$:
\begin{equation}
  \label{eq:18}
  \widetilde{\omega}_i^{j} =
  \begin{cases}
    \left[\omega_i^{j}\right] + (\nu_i^j-1)\tau^{-1}\left[\omega^{j-1}_{i}\right] &
    \text{if
      $\nu_i^j>1$}\\
    \left[\omega_i^{j}\right] & \text{if $\nu_i^j\leq 1$}
  \end{cases}
\end{equation}
Notice that, a priori it could happen that $\nu_i^0>1$ and the above definition would not
be correct. But this does not happen as the following lemma shows:
\begin{lemma}
  \label{prop:12}
  We have $\nu_i^0\leq 1$ for all $i=1,\dots,\mu_2$.
\end{lemma}
\begin{proof}
  The socle of the $0$-dimensional complete intersection defined by $(f)+I$ has degree
  $1+\mathbf{e}-\mathbf{p}$. Hence all the elements of degree greater than
  $1+\mathbf{e}-\mathbf{p}$ are contained in the image of the multiplication by $f$ and
  the lemma follows.
\end{proof}
For $\omega\in\omega_g$ let us set
\begin{displaymath}
  \lambda(\omega):=
  \begin{cases}
    1 & \text{if $\nu(\omega)>1$}\\
    \nu & \text{if $0\leq\nu(\omega)\leq 1$}\\
    0 & \text{if $\nu(\omega)<0$}
  \end{cases}
\end{displaymath}
and $\lambda(\omega_i^j):=\lambda_i^j$. In the next result we compute the spectral numbers
of $(\mathbf{G}_t,G_t)$ for $t=0$.
\begin{theorem}
  \label{prop:13}
  For any $\lambda\in\set{\lambda_i^j:1\leq i\leq\mu_2, 0\leq j \leq N_i}$ the classes of
  $\widetilde{\omega}_i^j$ for which $\lambda(\widetilde{\omega}_i^j)=\lambda$ induce a
  basis of the vector space $\text{gr}_\lambda^V(G_0/\tau^{-1}G_0)$. Hence the numbers
  $\lambda_i^j$ together with its multiplicities form the spectrum of
  $(\mathbf{G}_0,G_0)$.
\end{theorem}
\begin{proof}
  We only need to check that $\tau\partial_\tau+\lambda$ is nilpotent on
  $\text{gr}_\lambda^V\mathbf{G}_0$ for $\lambda\in[0,1]$. Notice first that by definition
  we have $f\omega_i^{N_i}\in I+(J)$. It follows from Lemma~\ref{prop:8} and that
  $\omega_i^{N_i}\in\Omega_{C_0}$. An straightforward calculation together with
  Cor.~\ref{prop:11} shows that
  \begin{align}
    \label{eq:19}
    \begin{split}
      \text{if $j<N_i$ then~} \tau\partial_\tau\widetilde{\omega}_i^j &=
      \begin{cases}
        \tau\widetilde{\omega}_i^{j+1} & \text{if~}\nu_i^j
        \leq 0 \\
        -\nu_i^j\widetilde{\omega}_i^j + \tau
        \widetilde{\omega}_i^{j+1} & \text{if~}0<\nu_i^j\leq 1 \\
        -\widetilde{\omega}_i^j + \tau \widetilde{\omega}_i^{j+1} & \text{if~}\nu_i^j>1,
      \end{cases}
      \\
      \text{and if $j=N_i$ }\tau\partial_\tau\widetilde{\omega}_i^{N_i} &=
      -\lambda_i^{N_i}\widetilde{\omega}_i^{N_i}.
    \end{split}
  \end{align}  
  We show the nilpotency of $\tau\partial_\tau+\lambda$ with some detail for the first
  case in (\ref{eq:19}) as the others are analogous. As $\nu_i^j\leq 0$ we have $j<N_i$ (see
  Remark~\ref{prop:9}). If $\nu_i^j<0$, then $\nu_i^{j+1}<1$ so that
  $\tau\partial_\tau\widetilde{\omega}_i^j\in V_{<0}(\mathbf{G})$. If $\nu_i^j=0$ then
  $\nu_i^{j+1}=1$ and we get
  \begin{equation}
    \label{eq:20}
    (\tau\partial_\tau)^2\widetilde{\omega}_i^j=
    \tau(\tau\partial_\tau+1)\widetilde{\omega}_i^{j+1}= 
    \begin{cases}
      \tau^2\widetilde{\omega}_i^{j+2} &
      \text{if $j+1<N_i$}\\
      0 & \text{if $j+1=N_i$}
    \end{cases}
  \end{equation}
  In both cases we have $(\tau\partial_\tau)^2\widetilde{\omega}_i^j\in
  V_{<0}(\mathbf{G})$.
\end{proof}
\begin{corollary}
  \label{prop:14}
  In the basis of $G_0$ induced by $\widetilde{\omega}_i^j$, the matrix of the action of
  $\partial_\tau$ takes the form
  \begin{equation}
    \label{eq:21}
    (A_0+A_\infty\tau^{-1})d\tau
  \end{equation}
  where $A_0$ and $A_\infty$ are constant matrices, and $A_\infty$ diagonal. In
  particular, $G_0$ extends to a bundle on $\pp^1$ with logarithmic connection on
  $\tau=0$.
\end{corollary}
In the case of plane curves, it turns out that the spectrum of the restriction $f_t$ of
$f$ to the fibre $C_t$ coincides with that of $C_0$. More precisely:
\begin{theorem}
  \label{prop:15}
  If $n=1$ then the classes of $\widetilde{\omega}_i^j$ induce a basis of
  $\text{gr}_\lambda^V(G_t/\tau^{-1}G_t)$ for any $t\in\cc$.
\end{theorem}
\begin{proof}
  The particularity of family of plane curves is that $(g)=I\subset J_g$. It follows as in
  the proof of the previous theorem that $\omega_i^{N_i}\in\Omega_{C_t}$. The proof of the
  theorem now follows almost verbatim, with the difference that we now have to use the
  full equation (\ref{eq:16}) in Lemma~\ref{prop:10}. For example, equation (\ref{eq:20})
  now becomes
  \begin{equation}
    \label{eq:22}
    (\tau\partial_\tau)^2\widetilde{\omega}_i^j=
    \tau(\tau\partial_\tau+1)\widetilde{\omega}_i^{j+1}= 
    \begin{cases}
      \tau^2\widetilde{\omega}_i^{j+2} &
      \text{if $j+1<N_i$}\\
      \tau t[\omega_{i,1}] + \tau^2t [\omega'_{i,1}] &
      \text{if $j+1=N_i$}
    \end{cases}
  \end{equation}
  which again in both cases belong to $V_{<0}(\mathbf{G}_t)$ as
  $\nu(\omega_{i,1}),\nu(\omega_{i,1}')\leq \nu_i^{N_i}-e_1=1-e_1$. The rest of cases
  are similarly adapted. Notice that for the elements $\widetilde{\omega}_i^{N_i}$ we
  might need to use Lemma~\ref{prop:10} say $K$ times, being $K=\min\set{k\geq 1:
    \nu_i^{N_i}-ke_1 <1}$ to ensure that
  $(\tau\nabla_\tau+\lambda_i^{N_i})^K\widetilde{\omega}_i^{N_i}\in
  V_{<\lambda_i^{N_i}}(\mathbf{G}_t)$.
\end{proof}
We can then use the results in \cite{DouSab03} to construct Frobenius manifolds on the
base space of the miniversal deformation of $f_t$ for $t\not=0$.
\begin{corollary}
  If $n=1$, the class of $\alpha$ in $G_t$ is a primitive form for any $t$. Hence for any
  $t\not=0$, the base space of the miniversal deformation of $f_t$ has the structure of a
  massive Frobenius manifold.
\end{corollary}
\begin{proof}
  Let $\omega_i^j=h_i^{(j)}\alpha$ be the basis of $G_0$ defined in \eqref{eq:17}. Then the
  unfolding $F=f+\sum_{i=1}^{\mu_2}\sum_{j=0}^{N_i} u_i^{(j)}h_i^{(j)}$ is miniversal. The
  connection with respect to the deformation parameters is given by
  \begin{equation}
    \label{eq:23}
    \widehat{\nabla}_{\partial_{u_i^{(j)}}}[\omega]
    =\left[\cvf[\omega]{u_i^{(j)}}\right]-\tau\left[\cvf[F]{u_i^{(j)}}\omega\right].
  \end{equation}
  As $\alpha=\omega_1^0$ we have
  \begin{align}
    \label{eq:24}
      \widehat{\nabla}_{\partial_{u_i^{(j)}}}[\alpha]=-\tau[\omega_i^j],~
      \widehat{\nabla}_{\partial_\tau}[\alpha]
      =[\omega_1^1]-\sum_{i=1}^{\mu_2}\sum_{j=0}^{N_i}u_i^{(j)}[\omega_i^j].
  \end{align}
  It follows that $\alpha$ is a primitive form. The existence of the metric follows from
  microlocal Poincar\'e duality (cf.~{\em loc.~cit.}). Finally, for a generic value of
  $u$, all the critical points of $F$ on $C_t\times\set{u}$ are Morse, hence the
  multiplication is generically semisimple.
\end{proof}
It is known that the metric is given by the sum of the residues at the critical
points. More precisely, if $(u_1,\dots,u_\mu)$ are parameters of the base space of the
miniversal deformation $F:(\cc^2\times B,0)\ra (\cc,0)$ and $dF=J_F\alpha$ denotes the
relative differential, then
\begin{equation}
  \label{eq:25}
  \gen{\cvf{u_i},\cvf{u_j}}_t
  =\int_{\partial C_t}\left.\left(\frac{\cvf[F]{u_i}\cvf[F]{u_j}}{J_F}\alpha\right)\right|_{C_t}
\end{equation}
where $\partial C_t$ is the boundary of an appropriate representative of the Milnor fibre
of $g$.
\begin{corollary}
  \label{prop:17}
  The formula \eqref{eq:25} for $t=0$ together with the multiplication defined by
  \eqref{eq:8} defines the structure of non-massive Frobenius manifold on $B$.
\end{corollary}
\begin{proof}
  We have $\gen{\cvf{u_i},\cvf{u_j}}_t\ra\gen{\cvf{u_i},\cvf{u_j}}_0$ when $t\ra 0$. The
  flatness of $\gen{-,-}_t$ implies that of $\gen{-,-}_0$ as it can be seen, for example,
  writing out explicit formulas for the curvature in terms of the Christoffel symbols. The
  existence of a potential can be translated into the flatness of the first structure
  connection (e.g.~\cite{Man99},~Th.~1.5). More precisely, for each $t$, let $^t\nabla$
  the Levi-Civita connection of $\gen{-,-}_t$. The first structure connection is defined
  as
  \begin{equation}
    \label{eq:26}
    ^t\overline{\nabla}_{z,\partial_{u_i}}\partial_{u_j}
    :=^t\nabla_{\partial_{u_i}}\partial_{u_j}+z\partial_{u_i}\star_t\partial_{u_j}
  \end{equation}
  It is of course closely related to the the Gauss-Manin connection
  $\widehat{\nabla}$. Notice that
  $\partial_{u_i}\star_t\partial_{u_j}\ra\partial_{u_i}\star_0\partial_{u_j}$ when $t\ra
  0$, and hence $^t\overline{\nabla}\ra{}^0\overline{\nabla}$. The result follows.
\end{proof}

\section{An example: linear functions on the $A_k$-singularity}
\label{sec:an-example}

Let us illustrate our construction with a worked-out example. We consider the curve $C_0$ defined by
$g(x,y)=x^k+y^2=0$, $k\geq 2$, and the function $f_0$ given by the restriction of $f(x,y)=x$ to $C_0$.

\vspace{\baselineskip}\noindent{\bf Miniversal deformation.} The classes of $1,\dots,x^{k-1}$ form a
$\cc$-basis of the Jacobian algebra $\oo_{C_0}/(2y)$ and hence a miniversal unfolding is given by $F+u_1 x^{k-1} +\dots+ u_{k-1}x + u_k$.

\vspace{\baselineskip}\noindent{\bf Spectrum}. For a homogeneous polynomial $h$ we have
\begin{displaymath}
  \nu(h)=\deg(h) - \frac{k-2}{2}
\end{displaymath}
According to Theorem~\ref{prop:13}, the spectrum of $f_t=f|_{C_t}$ is
\begin{align}
  \label{eq:32}
  \begin{split}
    &\set{\left(0,\frac{k}{2}\right),\left(1,\frac{k}{2}\right)}
    ~\text{if $k$ is even and,} \\
    &
    \set{\left(0,\frac{k-1}{2}\right),\left(\frac{1}{2},1\right),\left(1,\frac{k-1}{2}\right)}
    ~\text{if $k$ is odd.}
  \end{split}
\end{align}

\vspace{\baselineskip}\noindent{\bf Nilpotent Frobenius structure}. If we set
$F'=\cvf[F]{x}$, the multiplication table on $\Theta_{B,0}$
is given by the isomorphism
\begin{equation}
  \label{eq:33}
  \partial_{u_i} \xra{t'\mathscr{F}_0} x^{k-i} \in
  \pi_*\left(\frac{\oo}{(x^k+y^2,2yF')}\right)
\end{equation}
where $\oo$ denotes the sheaf of holomorphic functions on the variables
$x,y,u_1,\dots,u_k$ and $\pi:C_0\times (B,0) \ra (B,0)$ is the canonical projection. The
ideal $(x^k+y^2,2yF')$ defines in $C_0\times B$ a scheme with two components:
$W_1:=\set{0}\times B$ and the (reduced) variety $W_2$ defined by $F'=0$. As
$W_1$ already has multiplicity $k=\mu$, the $F$-manifold structure extends to
$B\setminus\pi(W_2)$ (notice that $0\not\in\pi(W_2)$). We see that this $F$-manifold
structure is purely nilpotent, in the sense that if $i\not=k$ (i.e., if $\partial_{u_i}$
is not the identity), we have $\partial_{u_i}\star\dots\star\partial_{u_i}=0$ where the
product occurs at most $k$ times. We remark that this is always the case if the function
$f_0$ is the restriction of a linear function as all the critical points are provided by
the singular curve.

The metric is also easy to describe, at least on $T_0 B$ (and hence on flat
coordinates). The generator of the socle of $\cc[x,y]/(x^k+y^2,2y)$ is $x^{k-1}$, so that
if we choose a residue form with $\Res(x^{k-1})=1$, the metric in the basis
$\partial_{u_i}|_0$ is simply given by the matrix with all entries equal $1$ in the
anti-diagonal, and $0$ everywhere else.


\begin{thebibliography}{10}

\bibitem{Bar00}
Serguei Barannikov, \emph{Semi-infinite {H}odge structures and mirror symmetry
  for projective spaces}, AG/0108148.

\bibitem{Bjork93}
Jan-Erik Bj{\"o}rk, \emph{Analytic {$\mathscr{D}$}-modules and applications},
  Mathematics and its Applications, vol. 247, Kluwer Academic Publishers Group,
  Dordrecht, 1993. \MR{MR1232191 (95f:32014)}

\bibitem{deGre06}
Ignacio de~Gregorio, \emph{Deformations of functions and {$F$}-manifolds}, To
  appear in Bull. London Math. Soc. (2006), math.AG/0503323.

\bibitem{DijVerVer91}
Robbert Dijkgraaf, Herman Verlinde, and Erik Verlinde, \emph{Notes on
  topological string theory and {$2$}{D} quantum gravity}, String theory and
  quantum gravity (Trieste, 1990), World Sci. Publishing, River Edge, NJ, 1991,
  pp.~91--156. \MR{93c:81202}

\bibitem{DouSab03}
A.~Douai and C.~Sabbah, \emph{Gauss-{M}anin systems, {B}rieskorn lattices and
  {F}robenius structures. {I}}, Proceedings of the International Conference in
  Honor of Fr\'ed\'eric Pham (Nice, 2002), vol.~53, 2003, pp.~1055--1116. \MR{2
  033 510}

\bibitem{Dub96}
Boris Dubrovin, \emph{Geometry of {$2$}{D} topological field theories},
  Integrable systems and quantum groups (Montecatini Terme, 1993), Lecture
  Notes in Math., vol. 1620, Springer, Berlin, 1996, pp.~120--348.
  \MR{97d:58038}

\bibitem{HerMan99}
C.~Hertling and Yu. Manin, \emph{Weak {F}robenius manifolds}, Internat. Math.
  Res. Notices (1999), no.~6, 277--286. \MR{2000j:53117}

\bibitem{Her02}
Claus Hertling, \emph{Frobenius manifolds and moduli spaces for singularities},
  Cambridge Tracts in Mathematics, vol. 151, Cambridge University Press,
  Cambridge, 2002. \MR{1 924 259}

\bibitem{Loo84}
E.~J.~N. Looijenga, \emph{Isolated singular points on complete intersections},
  London Mathematical Society Lecture Note Series, vol.~77, Cambridge
  University Press, Cambridge, 1984. \MR{86a:32021}

\bibitem{Man99}
Yuri~I. Manin, \emph{Frobenius manifolds, quantum cohomology, and moduli
  spaces}, American Mathematical Society Colloquium Publications, vol.~47,
  American Mathematical Society, Providence, RI, 1999. \MR{MR1702284
  (2001g:53156)}

\bibitem{Mann05}
{\'E}.~Mann, \emph{Cohomologie quantique orbifolde des espaces projectifs {\`a}
  poids}, Ph.D. thesis, Universit{\'e} Louis Pasteur, 2005.

\bibitem{MonStr01}
David Mond and Duco van Straten, \emph{Milnor number equals {T}jurina number
  for functions on space curves}, J. London Math. Soc. (2) \textbf{63} (2001),
  no.~1, 177--187. \MR{2002e:32040}

\bibitem{Sab98}
Claude Sabbah, \emph{Frobenius manifolds: isomonodromic deformations and
  infinitesimal period mappings}, Exposition. Math. \textbf{16} (1998), no.~1,
  1--57. \MR{MR1617534 (99k:32031)}

\bibitem{KSai81}
Kyoji Saito, \emph{Primitive forms for a universal unfolding of a function with
  an isolated critical point}, J. Fac. Sci. Univ. Tokyo Sect. IA Math.
  \textbf{28} (1981), no.~3, 775--792 (1982). \MR{MR656053 (84k:32031)}

\bibitem{KSai83}
\bysame, \emph{Period mapping associated to a primitive form}, Publ. Res. Inst.
  Math. Sci. \textbf{19} (1983), no.~3, 1231--1264. \MR{85h:32034}

\bibitem{MSai89}
Morihiko Saito, \emph{On the structure of {B}rieskorn lattice}, Ann. Inst.
  Fourier (Grenoble) \textbf{39} (1989), no.~1, 27--72. \MR{91i:32035}

\bibitem{CoaCorLeeTse06}
Yuan-Pin Lee Hsian-Hua~Tseng Tom~Coates, Alessio~Corti, \emph{The quantum
  orbifold cohomology of weighted projective space}, math.AG/0608481.

\end{thebibliography}
\def\cprime{$'$}
\providecommand{\bysame}{\leavevmode\hbox to3em{\hrulefill}\thinspace}
\providecommand{\MR}{\relax\ifhmode\unskip\space\fi MR }
\providecommand{\MRhref}[2]{%
  \href{http://www.ams.org/mathscinet-getitem?mr=#1}{#2}
}
\providecommand{\href}[2]{#2}

\end{document}